\documentclass{article}%
\usepackage{amsmath}
\usepackage{amssymb}
\usepackage{hyperref}
\usepackage{amsfonts}
\usepackage{graphicx}%
\begin{document}

\title{Bernoulli and Euler Partitions}
\author{Thomas Curtright$^{\S }$ and Christophe Vignat$^{\circledcirc}$
\and \ \ \ \ \ \ \ {\footnotesize curtright@miami.edu
\ \ \ \ \ \ \ \ \ \ \ \ \ \ \ \ cvignat@tulane.edu \ \ \ \ \ \ \ \ \ \ \ }%
\medskip\\$^{\S }$Department of Physics, University of Miami, Coral Gables, FL 33124\\$^{\circledcirc}$Department of Mathematics, Tulane University, New Orleans, LA 70118}
\date{\vspace{-0.25in}}
\maketitle

\begin{abstract}
Exact rational partitions are presented for Bernoulli and Euler numbers as
novel sums involving Faulhaber and Sali\'{e} coefficients.

\end{abstract}

Some interesting well-ordered and strictly positive partitions\ \cite{TLC}
for\footnote{For comparison with the Bernoulli numbers, we choose to omit the
$2^{2n}$ factor used in the conventional definition of Euler numbers in terms
of Euler polynomials, namely, $E_{2n}\equiv2^{2n}E_{2n}\left(  x=1/2\right)
$.} $\left\vert E_{2n}\left(  1/2\right)  \right\vert =\left(  -1\right)
^{n}E_{2n}\left(  x=1/2\right)  $ and $\left\vert B_{2n}\right\vert =\left(
-1\right)  ^{n+1}B_{2n}\left(  x=0\right)  $ may be obtained by writing matrix
equations,
\begin{equation}
\overrightarrow{1}=M_{E}\cdot\overrightarrow{E}=M_{B}\cdot\overrightarrow{B}
\label{StartYerEngines!}%
\end{equation}
where $\overrightarrow{E}$ and $\overrightarrow{B}$\ are columns of these even
index Euler and Bernoulli polynomial values, $\overrightarrow{1}$\ is a column
of $1$s, and the lower triangular matrices $M_{E,B}$ have entries
\begin{equation}
\left(  M_{E}\right)  _{m,n}=\left(  -1\right)  ^{m-n}4^{m}\binom{m}{2\left(
m-n\right)  }\ ,\ \ \ \left(  M_{B}\right)  _{m,n}=2\left(  -1\right)
^{m-n}\binom{2n-1}{m}\binom{2m+1}{2n}%
\end{equation}
Computing the lower triangular inverses for these matrices then decomposes
each of the above Euler and Bernoulli polynomial values as a finite sum of
positive rational numbers (\textquotedblleft partitions\textquotedblright). \
\begin{equation}
\overrightarrow{E}=M_{E}^{-1}\cdot\overrightarrow{1}%
\ ,\ \ \ \overrightarrow{B}=M_{B}^{-1}\cdot\overrightarrow{1}
\label{E&BPartitions}%
\end{equation}
To illustrate this, the first six rows and columns of the triangular matrices
$M_{B,E}$ and $M_{B,E}^{-1}$ are given by%
\[
\left(  M_{B}\right)  \boldsymbol{=}\left(
\begin{array}
[c]{cccccc}%
6 & 0 & 0 & 0 & 0 & 0\\
0 & 30 & 0 & 0 & 0 & 0\\
0 & -70 & 140 & 0 & 0 & 0\\
0 & 0 & -840 & 630 & 0 & 0\\
0 & 0 & 924 & -6930 & 2772 & 0\\
0 & 0 & 0 & 18\,018 & -48\,048 & 12\,012
\end{array}
\right)  \ ,\ \ \left(  M_{B}^{-1}\right)  =\left(
\begin{array}
[c]{cccccc}%
\frac{1}{6}\smallskip & 0 & 0 & 0 & 0 & 0\\
0 & \frac{1}{30}\smallskip & 0 & 0 & 0 & 0\\
0 & \frac{1}{60}\smallskip & \frac{1}{140} & 0 & 0 & 0\\
0 & \frac{1}{45}\smallskip & \frac{1}{105} & \frac{1}{630} & 0 & 0\\
0 & \frac{1}{20}\smallskip & \frac{3}{140} & \frac{1}{252} & \frac{1}{2772} &
0\\
0 & \frac{1}{6}\smallskip & \frac{1}{14} & \frac{17}{1260} & \frac{1}{693} &
\frac{1}{12\,012}%
\end{array}
\right)
\]%
\[
\left(  M_{E}\right)  \boldsymbol{=}\left(
\begin{array}
[c]{cccccc}%
4 & 0 & 0 & 0 & 0 & 0\\
-16 & 16 & 0 & 0 & 0 & 0\\
0 & -192 & 64 & 0 & 0 & 0\\
0 & 256 & -1536 & 256 & 0 & 0\\
0 & 0 & 5120 & -10\,240 & 1024 & 0\\
0 & 0 & -4096 & 61\,440 & -61\,440 & 4096
\end{array}
\right)  \ ,\ \ \left(  M_{E}^{-1}\right)  \boldsymbol{=}\left(
\begin{array}
[c]{cccccc}%
\frac{1}{4}\smallskip & 0 & 0 & 0 & 0 & 0\\
\frac{1}{4}\smallskip & \frac{1}{16} & 0 & 0 & 0 & 0\\
\frac{3}{4}\smallskip & \frac{3}{16} & \frac{1}{64} & 0 & 0 & 0\\
\frac{17}{4}\smallskip & \frac{17}{16} & \frac{3}{32} & \frac{1}{256} & 0 &
0\\
\frac{155}{4}\smallskip & \frac{155}{16} & \frac{55}{64} & \frac{5}{128} &
\frac{1}{1024} & 0\\
\frac{2073}{4}\smallskip & \frac{2073}{16} & \frac{23}{2} & \frac{135}{256} &
\frac{15}{1024} & \frac{1}{4096}%
\end{array}
\right)
\]
hence the partitions of $\left\vert E_{2n}\left(  1/2\right)  \right\vert $
and $\left\vert B_{2n}\right\vert $ for $n=1$ to $6$. \ For example,
$\left\vert E_{12}\left(  1/2\right)  \right\vert =\frac{2702\,765}%
{4096}=\frac{2073}{4}+\frac{2073}{16}+\frac{23}{2}+\frac{135}{256}+\frac
{15}{1024}+\frac{1}{4096}$ and $\left\vert B_{12}\right\vert =\frac{691}%
{2730}=\frac{1}{6}+\frac{1}{14}+\frac{17}{1260}+\frac{1}{693}+\frac
{1}{12\,012}$.\newpage

Explicit general results for the partitions implicit in (\ref{E&BPartitions})
are presented in the following, sans proofs. \ Detailed proofs are available
and will be given later in a longer exposition of this work \cite{CandV}%
.\bigskip

For the $6\times6$ matrices shown above, various properties of the entries in
$M_{E,B}^{-1}$ can be easily checked. \ For example, the diagonals of
$M_{E,B}^{-1}$ are given by $1/4^{n}$ and $\frac{\left(  n!\right)  ^{2}%
}{\left(  2n+1\right)  !}$, respectively, where $n$ is the row number. \ The
first sub-diagonals for $n\geq2$ are respectively $2n\left(  n-1\right)
/4^{n}$\ and $\frac{1}{6}\left(  n-2\right)  \frac{n!\left(  n-1\right)
!}{\left(  2n-1\right)  !}$, the second sub-diagonals for $n\geq3$ are
respectively $\frac{2}{3}n\left(  n-1\right)  \left(  n-2\right)  \left(
5n-3\right)  /4^{n}$\ and $\frac{7}{360}\left(  n-3\right)  \left(  n-\frac
{8}{7}\right)  \frac{n!\left(  n-2\right)  !}{\left(  2n-3\right)  !}$,
etc.\bigskip

In general, the even index Euler values decompose as the finite partitions%
\begin{equation}
\left\vert E_{2n}\left(  1/2\right)  \right\vert =\left(  -1\right)
^{n}E_{2n}\left(  x=1/2\right)  =\sum_{k=1}^{n}\epsilon\left(  n,k\right)
\text{ \ for \ }n\geq1 \label{EulerPartition}%
\end{equation}
where $\left(  M_{E}^{-1}\right)  _{n,k}\equiv\epsilon\left(  n,k\right)
\geq0$ for all integer $n\geq k\geq0$. \ The $\epsilon$s are rescaled moduli
of \textquotedblleft Sali\'{e} coefficients\textquotedblright%
\ \cite{Salie,GandV}, $s\left(  n,k\right)  =\left(  -1\right)  ^{n-k}%
\left\vert s\left(  n,k\right)  \right\vert $ with $\epsilon\left(
n,k\right)  =\left\vert s\left(  n,k\right)  \right\vert /4^{k}$, where the
$s\left(  n,k\right)  $ may be obtained from a bivariate generating function
\cite{GandV}.%
\begin{equation}
\sum_{k,n\geq0}\frac{x^{2n}}{\left(  2n\right)  !}~s\left(  n,k\right)
y^{k}=\frac{\cosh\left(  \frac{1}{2}~x\sqrt{1+4y}\right)  }{\cosh\left(
\frac{1}{2}~x\right)  }%
\end{equation}
The Sali\'{e} coefficients appear in alternating sums of integer powers of the
first $m>0$ integers, when re-expressed as powers of $m\left(  m+1\right)  $.
\ For example \cite{GandV}, if $n>0$,%
\begin{equation}
\sum_{l=1}^{m}\left(  -1\right)  ^{m-l}l^{2n}=\frac{1}{2}\sum_{k=0}%
^{n}s\left(  n,k\right)  \left(  m\left(  m+1\right)  \right)  ^{k}%
\end{equation}

For comparison, the even index Bernoulli numbers decompose as the finite
partitions\footnote{For $B_{2n}$, the $n=1$ case is not partitioned here, like
the $n=2$ case, as is evident in the $6\times6$ matrices given previously.
\ Namely, $B_{2}=1/6$ and $-B_{4}=1/30$. \ This feature is shared by $E_{0}=1$
and $E_{2}\left(  x=1/2\right)  =1/4$. \ This, as well as the different
summation lower limits in (\ref{EulerPartition}) and (\ref{BernoulliPartition}%
), suggests the indices should be shifted when individual Bernoulli and Euler
partitions are compared.}
\begin{equation}
\left\vert B_{2n}\right\vert =\left(  -1\right)  ^{n+1}B_{2n}\left(
x=0\right)  =\sum_{k=2}^{n}b\left(  n,k\right)  \text{ \ for \ }n\geq2
\label{BernoulliPartition}%
\end{equation}
where $\left(  M_{B}^{-1}\right)  _{n,k}\equiv b\left(  n,k\right)  \geq0$ for
all integer $n\geq k\geq1$.\ The $b$s are moduli of rescaled and shifted
\textquotedblleft Faulhaber coefficients\textquotedblright\ where the latter
are\ designated as $A_{m}^{\left(  n\right)  }$ by Knuth \cite{Knuth} or as
$f\left(  n,m\right)  \equiv A_{n-m}^{\left(  n+1\right)  }/\left(
n+1\right)  $ by Gessel and Viennot \cite{GandV}. \ Explicitly, for $n\geq
k\geq2$,%
\begin{equation}
b\left(  n,k\right)  =\left(  -1\right)  ^{n-k}n~f\left(  n-1,k-1\right)
~\frac{\left(  k!\right)  ^{2}}{\left(  2k+1\right)  !}=\left(  -1\right)
^{n-k}A_{n-k}^{\left(  n\right)  }~\frac{\left(  k!\right)  ^{2}}{\left(
2k+1\right)  !} \label{bFaul}%
\end{equation}
The Faulhaber coefficients may be obtained from bivariate generating
functions, either as \cite{Knuth}
\begin{equation}
\sum_{n\geq1}^{\infty}\sum_{k=0}^{n-1}\frac{x^{2n}}{\left(  2n\right)  !}%
A_{k}^{\left(  n\right)  }y^{k}=\left(  \frac{x\sqrt{y}}{2}\right)
\frac{\cosh\left(  \frac{1}{2}x\sqrt{y+4}\right)  -\cosh\left(  \frac{1}%
{2}x\sqrt{y}\right)  }{\sinh\left(  \frac{1}{2}x\sqrt{y}\right)  }%
\end{equation}
or else as \cite{GandV}%
\begin{equation}
\sum_{n,k\geq0}\frac{x^{2n}}{\left(  2n\right)  !}~f\left(  n,k\right)
y^{k}=\frac{d}{dx}\left(  \frac{\cosh\left(  \frac{1}{2}x\sqrt{1+4y}\right)
-\cosh\left(  \frac{1}{2}x\right)  }{y\sinh\left(  \frac{1}{2}x\right)
}\right)
\end{equation}
That is to say,%
\begin{equation}
\sum_{n,k\geq0}\frac{x^{2n}}{\left(  2n\right)  !}~f\left(  n,k\right)
y^{k}=\frac{\left(  \sqrt{1+4y}-1\right)  }{2y}\frac{\sinh\left(  \frac{1}%
{2}x\sqrt{1+4y}\right)  }{\sinh\left(  \frac{1}{2}x\right)  }+\frac
{1-\cosh\left(  \frac{1}{2}x\left(  \sqrt{1+4y}-1\right)  \right)  }{2y\left(
\sinh\left(  \frac{1}{2}x\right)  \right)  ^{2}}%
\end{equation}
The $f\left(  n,k\right)  $ appear in sums of integer powers of the first
$m>0$ integers, when re-expressed as powers of $m\left(  m+1\right)  $. For
example,
\begin{equation}
\sum_{l=1}^{m}l^{2n+1}=\frac{1}{2}\sum_{k=0}^{n}f\left(  n,k\right)  \left(
m\left(  m+1\right)  \right)  ^{k+1}%
\end{equation}
a feature first noted by Faulhaber in the 17$^{th}$ century \cite{Faulhaber}.

Alternatively, the $b$s are conveniently evaluated by a computer upon
expressing the Faulhaber coefficients in terms of $\Gamma$ functions and
Riemann $\zeta$s. \ For $n\geq k\geq1$,%
\begin{equation}
A_{n-k}^{\left(  n+1\right)  }=\left(  n+1\right)  f\left(  n,k\right)
=\frac{2\left(  -1\right)  ^{n-k}\Gamma\left(  2n+3\right)  }{\Gamma\left(
k+2\right)  }\sum_{j=0}^{\left\lfloor \left(  k-1\right)  /2\right\rfloor
}\frac{\left(  -1\right)  ^{j}\Gamma\left(  2k-2j\right)  }{\Gamma\left(
k-2j\right)  \Gamma\left(  2j+2\right)  }\frac{\zeta\left(  2n-2j\right)
}{\left(  2\pi\right)  ^{2n-2j}} \label{f as Gamma Zeta}%
\end{equation}
which leads, for $n\geq k\geq2$, to the result%
\begin{equation}
b\left(  n,k\right)  =\frac{4\Gamma\left(  2n+1\right)  \Gamma\left(
k+2\right)  }{\Gamma\left(  2k+3\right)  }\sum_{j=0}^{\left\lfloor
k/2\right\rfloor -1}\frac{\left(  -1\right)  ^{j}\Gamma\left(  2k-2-2j\right)
}{\Gamma\left(  k-1-2j\right)  \Gamma\left(  2j+2\right)  }\frac{\zeta\left(
2n-2j-2\right)  }{\left(  2\pi\right)  ^{2n-2j-2}} \label{b as Gamma Zeta}%
\end{equation}
This form for $b\left(  n,k\right)  $ immediately yields the asymptotic
behavior as $n\rightarrow\infty$ for fixed $k$ (previously conjectured in
\cite{TLC}) as a simple consequence of $\zeta\left(  m\right)
\underset{m\rightarrow\infty}{\sim}1$. \ 

Since $\left\vert B_{2m}\right\vert =2\left(  2m\right)  !\zeta\left(
2m\right)  /\left(  4\pi^{2}\right)  ^{m}$, as a consequence of
(\ref{b as Gamma Zeta}) the Bernoulli partitions may be re-expressed as a
recursion relation giving $\zeta\left(  n\right)  $ values for even $n$. For
$m\geq1$,%
\begin{equation}
\zeta\left(  2m+2\right)  =\sum_{n=1}^{m}\sum_{k=1}^{\left\lfloor \left(
n+1\right)  /2\right\rfloor }\frac{\sqrt{4\pi}\left(  -1\right)  ^{k-1}%
\pi^{2k}\Gamma\left(  2n-2k+2\right)  \zeta\left(  2m+2-2k\right)  }%
{4^{n-k+2}\Gamma\left(  n+\frac{5}{2}\right)  \Gamma\left(  2k\right)
\Gamma\left(  n+2-2k\right)  }\label{ZetaRecursion}%
\end{equation}
This recursion is not easily found in the literature\footnote{All the $\zeta$s
in (\ref{ZetaRecursion}) may be supplanted by even indexed $B$s and the
partition of the Bernoulli numbers is thereby recast as a recursion relation
expressing $B_{2m}$ in terms of lesser indexed $B$s. \ The result is
\cite{CandV}%
\[
B_{2m}=\frac{-1}{\left(  m+1\right)  \left(  2m+1\right)  }\sum
_{j=m-\left\lfloor m/2\right\rfloor }^{m-1}\binom{m+1}{2m-2j+1}\left(
1+2j\right)  B_{2j}\text{ \ \ for \ \ }m\geq2
\]
This recursion relation seems to be new, but it can be expressed as a linear
combination of some novel relations obtained by Apostol (see Eqn(38) in
\cite{Apostol}). \ Suffice it to say here, this new recursion provides a
direct route that leads to the solution of the Bernoulli relation in
Eqn(\ref{StartYerEngines!}). \ That said, using Riemann $\zeta$s instead of
Bernoulli numbers in (\ref{ZetaRecursion}) not only makes the asymptotic
behavior of $b\left(  n,k\right)  $ transparent, it also leads naturally to
the infinite sum and integral relations to follow, namely,
(\ref{sums and integrals}).}. \ Also, because $b\left(  n,k\right)  \geq0$ for
all integer $n\geq k\geq1$, when summed over $k$ each $n$-summand in
(\ref{ZetaRecursion}) is positive.

The expressions for $f\left(  n,k\right)  $\ or $b\left(  n,k\right)  $\ in
(\ref{f as Gamma Zeta}) or (\ref{b as Gamma Zeta}) can be recast as infinite
sums or integrals, with modified summands and integrands, upon using the
corresponding well-known expressions for $\zeta$,%
\begin{equation}
\zeta\left(  2m\right)  =\sum_{l=1}^{\infty}\frac{1}{l^{2m}}=\frac{1}%
{\Gamma\left(  2m\right)  }\int_{0}^{\infty}\frac{s^{2m-1}}{e^{s}-1}~ds
\end{equation}
The corresponding infinite sum and integral forms of $b\left(  n,k\right)  $
for $2\leq k\leq n$ are given by%
\begin{equation}
b\left(  n,k\right)  =\sum_{l=1}^{\infty}p_{1}\left(  n,k,l\right)  =\int%
_{0}^{\infty}\frac{1}{e^{2\pi t}-1}~p_{2}\left(  n,k,t\right)  dt
\label{sums and integrals}%
\end{equation}
where $p_{1}\left(  n,k,1/l\right)  $ is a hypergeometric polynomial in $1/l$
of order $2n-2$,%
\begin{equation}
p_{1}\left(  n,k,1/l\right)  =\frac{4}{\left(  2\pi l\right)  ^{2n-2}}%
\frac{\Gamma\left(  2n+1\right)  \Gamma\left(  k+2\right)  }{\Gamma\left(
2k+3\right)  }\sum_{j=0}^{\left\lfloor k/2\right\rfloor -1}\frac{\left(
-1\right)  ^{j}\Gamma\left(  2k-2-2j\right)  }{\Gamma\left(  k-1-2j\right)
\Gamma\left(  2j+2\right)  }\left(  2\pi l\right)  ^{2j} \label{p1}%
\end{equation}
and where $p_{2}\left(  n,k,t\right)  ~$is a hypergeometric polynomial in $t$
of order $2n-3$,%
\begin{equation}
p_{2}\left(  n,k,t\right)  =4t^{2n-3}\frac{\Gamma\left(  2n+1\right)
\Gamma\left(  k+2\right)  }{\Gamma\left(  2k+3\right)  }\sum_{j=0}%
^{\left\lfloor k/2\right\rfloor -1}\frac{\left(  -1\right)  ^{j}\Gamma\left(
2k-2-2j\right)  }{\Gamma\left(  k-1-2j\right)  \Gamma\left(  2j+2\right)
\Gamma\left(  2n-2j-2\right)  }\left(  \frac{1}{t}\right)  ^{2j} \label{p2}%
\end{equation}
The lowest power of $1/l$ in $p_{1}\left(  n,k,t\right)  $ is $\left(
1/l\right)  ^{2n-2\left\lfloor k/2\right\rfloor }$, hence $p_{1}\left(
n,k,1/l\right)  =\left(  1/l\right)  ^{2n-2\left\lfloor k/2\right\rfloor
}q_{1}\left(  n,k,1/l\right)  $ where $q_{1}\left(  n,k,t\right)  $ for $2\leq
k\leq n$ is a polynomial in $1/l$ of order $2\left\lfloor \frac{1}%
{2}k\right\rfloor -2$. \ The lowest power of $t$ in $p_{2}\left(
n,k,t\right)  $ is $t^{2n-1-2\left\lfloor k/2\right\rfloor }$, hence
$p_{2}\left(  n,k,t\right)  =t^{2n-1-2\left\lfloor k/2\right\rfloor }%
q_{2}\left(  n,k,t\right)  $ where $q_{2}\left(  n,k,t\right)  $ for $2\leq
k\leq n$ is a polynomial in $t$ of order $2\left\lfloor \frac{1}%
{2}k\right\rfloor -2$. \ 

Either Maple or Mathematica can readily express these $p_{1,2}$ and $q_{1,2}$
polynomials in terms of generalized hypergeometric functions whose sums and
weighted integrals therefore reduce to Faulhaber coefficients as given by
(\ref{sums and integrals}) upon using (\ref{bFaul}). \ So expressed, those
summation and integration results are elusive in the literature, if not
altogether absent.

Next, consider the sense in which the terms in the partitions are
well-ordered. \ It is not difficult to show for the Bernoulli partitions%
\begin{equation}
b\left(  n,k\right)  >2~b\left(  n,k+1\right)  \label{order}%
\end{equation}
for $n>2$ and $2\leq k\leq n-1$, and therefore $b\left(  n,k\right)  >0$ since
$b\left(  n,n\right)  =\left(  n!\right)  ^{2}/\left(  2n+1\right)  !>0$.
\ Other consequences of (\ref{order}) follow immediately. \ For example, again
for $n>2$ and $2\leq k\leq n-1$,%
\begin{equation}
b\left(  n,k\right)  >b\left(  n,k+1\right)  +2b\left(  n,k+2\right)
>\cdots>\left(  \sum_{l=k+1}^{n}b\left(  n,l\right)  \right)  +b\left(
n,n\right)  >\sum_{l=k+1}^{n}b\left(  n,l\right)  \label{tail}%
\end{equation}
The same inequalities are obeyed by the $\epsilon\left(  n,k\right)  $ in the
Euler partitions, for $n>1$ and $1\leq k\leq n-1$. \ 

Finally, the full set of Bernoulli partitions and inequalities may be recast
in terms of the ratios $r\left(  n,k\right)  \equiv b\left(  n,k\right)
/\left\vert B_{2n}\right\vert $, for $2\leq k\leq n$, such that $\sum
_{k=2}^{n}r\left(  n,k\right)  =1$. \ So, from (\ref{tail}),%
\begin{equation}
r\left(  n,k\right)  >\sum_{l=k+1}^{n}r\left(  n,l\right)  \text{ \ \ for
}n>2\text{ and }2\leq k\leq n-1.\
\end{equation}
That is to say, for each $n$ the $r\left(  n,k\right)  $ provide a partition
of unity as a sum of $n-1$ well-ordered, strictly positive, rational numbers,
thereby allowing an interpretation as a monotonic probability distribution.
\ Corresponding statements apply to the Euler case. \ 

To close, there are two obvious questions to ask. \ (1) What are the
combinatoric interpretations of these partitions? \ This is an open question.
\ (2) What are the applications of the results presented here? \ This too is
an open question for the Euler partitions, but for the Bernoulli partitions
there is at least one application known: \ Scattering by a scale invariant
potential \cite{CV,TLCagain}. \ In fact, this particular application
stimulated the work and results presented here.

\end{document}